\documentclass[12pt]{article}
\usepackage{amssymb, amsfonts, amssymb, amsthm, amscd}

\newtheorem{theorem}{{\bf Theorem}}[section]
\newtheorem{lemma}{{\bf Lemma}}[section]
\newtheorem{proposition}{{\bf Proposition}}[section]
\newtheorem{definition}{{\bf Definition}}[section]

\newtheorem{remark}{{\bf Remark}}[section]

\begin{document}

\title{Construction and Uniqueness for reflected BSDE under linear increasing
condition}
\author{G. Jia$^{a}$ and\ \ \
Mingyu XU $^{b}$\footnote{Corresponding author, Email:
xvmingyu@gmail.com}
\\
{\small $^a$School of Mathematics and System Science, Shandong University, }%
\\
{\small 250100, Jinan, China}\\
{\small $^b$Department of Financial Mathematics and Control science,
School
of Mathematical Science,} \\
{\small Fudan University, Shanghai, 200433, China.}}
\date{This version is done on April 2006.}
\maketitle

\textbf{Abstract.}In this paper, we study the uniqueness of the
solution of reflected BSDE with one or two barriers, under
continuous and linear increasing condition of generator $g$. Before
that we study the construction of solution of of reflected BSDE with
one or two barriers.

\section{Introduction}

El Karoui, Kapoudjian, Pardoux, Peng and Quenez introduced in 1997
the notion of reflected BSDE (RBSDE in short) on one lower barrier
\cite{EKPPQ}: the solution is forced to remain above a continuous
process, which is considered as the lower barrier. More precisely, a
solution for such equation associated to a coefficient $f$, a
terminal value $\xi $, a continuous barrier $L$, is a triple
$(Y_{t},Z_{t},K_{t})_{0\leq t\leq T}$ of adapted processes valued on
$\mathbb{R}^{1+d+1}$, which satisfies a square integrability
condition,
\[
Y_{t}=\xi
+\int_{t}^{T}f(s,Y_{s},Z_{s})ds+K_{T}-K_{t}-\int_{t}^{T}Z_{s}dB_{s},0\leq
t\leq T\mbox{, a.s.,}
\]
and $Y_{t}\geq L_{t}$, $0\leq t\leq T$, a.s.. Furthermore, the process $%
(K_{t})_{0\leq t\leq T}$ is non decreasing, continuous, and the role of $%
K_{t}$ is to push upward the state process in a minimal way, to keep
it above $L$. In this sense it satisfies
$\int_{0}^{T}(Y_{s}-L_{s})dK_{s}=0$.
They proved existence and uniqueness of a solution when $f$ is Lipschitz in $%
(y,z)$ uniformly in $(t,\omega )$. Then Matoussi (1997, \cite{M97})
considered the case $f$ continuous and at most linear growth in $y$,
$z$ and proved the existence of a maximal and a minimal solution.

Cvitanic and Karatzas (1995) studied the backward stochastic
differential equation with two barriers. A solution to such
equation associated to a terminal condition $\xi $, a coefficient $%
f(t,\omega ,y,z)$ and two barriers $L$ and $U$, is a triple
$(Y,Z,K)$ of adapted processes, valued in $\mathbf{R}^{1+d+1}$,
which satisfies $\;$
\[
Y_{t}=\xi
+\int_{t}^{T}f(s,Y_{s},Z_{s})ds+K_{T}^{+}-K_{t}^{+}-(K_{T}^{-}-K_{t}^{-})-%
\int_{t}^{T}Z_{s}dB_{s},\;0\leq t\leq T\mbox{ a.s}.
\]
$L_{t}\leq Y_{t}\leq U_{t}$, $0\leq t\leq T$ and $%
\int_{0}^{T}(Y_{s}-L_{s})dK_{s}^{+}=\int_{0}^{T}(Y_{s}-U_{s})dK_{s}^{-}=0,$
a.s. In this case, a solution $Y$ has remain between the lower
boundary $L$ and upper boundary $U$, almost surely. This is achieved
by the cumulative action of two continuous, increasing reflecting
processes $K^{\pm }$, which act in a minimal way when $Y$ attempts
to cross barriers.In this paper, authors proves the existence and
uniqueness of the solution, under certain condition of $\xi$, $L$
and $U$, and Lipschitz condition of generator $g$. In 1997,
Hamadene, Lepeltier and Matoussi extended the existence result to
the case when generator $g$ is a continuous function.

In this paper we study the uniqueness of the solution of reflected
BSDE with one or two barriers, under continuous and linear
increasing condition of generator $g$. Before that we study the
construction of solution of of reflected BSDE with one or two
barriers.

The paper is organized as following: in following section, we
present the basic assumption and notation of reflected BSDE with one
or two barriers, then we give existence result of reflected BSDE
with one or two barriers under linear increasing condition of $g$ in
section 3. In section 4 and 5, we study  the construction of
solution of of reflected BSDE with one or two barriers under linear
increasing condition of $g$. In section 6, we give a sufficient
result of uniqueness result of reflected BSDE with one or two
barriers. In section 7, we consider continuous depending on $c$ of
coefficient $g+c$.

\section{Assumptions and Notations}

Let $(\Omega ,\mathcal{F},P)$ be a complete probability space, and $%
(B_{t})_{0\leq t\leq T}=(B_{t}^{1},B_{t}^{2},\cdots ,B_{t}^{d})_{0\leq t\leq
T}^{\prime }$ be a $d$-dimensional Brownian motion defined on a finite
interval $[0,T]$, $0<T<+\infty $. Denote by $\{\mathcal{F}_{t};0\leq t\leq
T\}$ the standard filtration generated by the Brownian motion $B$, i.e. $%
\mathcal{F}_{t}$ is the completion of
\[
\mathcal{F}_{t}=\sigma \{B_{s};0\leq s\leq t\},
\]
with respect to $(\mathcal{F},P)$. We denote by $\mathcal{P}$ the $\sigma $%
-algebra of predictable sets on $[0,T]\times \Omega $.

We shall need the following spaces:

\[
\begin{array}{ll}
\mathbf{L}^{2}(\mathcal{F}_{t})= & \{\eta
:\mathcal{F}_{t}\mbox{-measurable
random real-valued variable, s.t. }E(|\eta |^{2})<+\infty \}, \\
\mathcal{H}_{n}^{2}(0,T)= & \{(\psi _{t})_{0\leq t\leq
T}:\mbox{predictable process valued in }\mathbb{R}^{n}\mbox{, s.t.
}E\int_{0}^{T}\left| \psi
(t)\right| ^{2}dt<+\infty \}, \\
\mathcal{S}^{2}(0,T)= & \{(\psi _{t})_{0\leq t\leq
T}:\mbox{progressively
measurable real-valued continuous process,} \\
& \mbox{s.t. }E(\sup_{0\leq t\leq T}\left| \psi (t)\right|
^{2})<+\infty \},
\\
\mathcal{A}^{2}(0,T)= & \{(K_{t})_{0\leq t\leq T}\in
\mathcal{S}^{2}(0,T):\ \mbox{increasing process, s.t. }K(0)=0\mbox{,
}E(K(T)^{2})<+\infty \}.
\end{array}
\]
To simplify the symbol, we use $\mathcal{H}^{2}(0,T)$ instead of $\mathcal{H}%
_{n}^{2}(0,T)$, when $n=1$.

To consider a reflected BSDE with one or two barriers, we are first given
\textbf{(i)} a terminal condition $\xi \in \mathbf{L}^{2}(\mathcal{F}_{T})$,
then \textbf{(ii)} a function
\[
g(t,\omega ,y,z):[0,T]\times \Omega \times \mathbf{R\times R}^{d}\rightarrow
\mathbf{R,}
\]
which is the coefficient of the reflected BSDE, and satisfies the followings:

(H1) for $(y,z)\in \mathbf{R\times R}^{d}$, $g(\cdot ,t,y,z)\in \mathcal{H}%
^{2}(0,T)$;

(H2) there exists a constant $\beta >0$, such that for $(t,\omega )\in
[0,T]\times \Omega $,
\begin{equation}
\left| g(t,y,z)\right| \leq \beta (1+\left| y\right| +\left| z\right| );
\label{linear}
\end{equation}

(H3) for $(t,\omega )\in [0,T]\times \Omega $, $g(t,\omega ,\cdot ,\cdot )$
is continuous.

And we need \textbf{(iii)} barriers $L$ and $U$, which are a progressively
measurable real-valued continuous process, such that
\begin{equation}
E[\sup_{0\leq t\leq T}(L_{t}^{+})^{2}+\sup_{0\leq t\leq
T}(U_{t}^{-})^{2}]<\infty \mbox{, and }L_{T}\leq \xi \leq U_{T}.
\label{barrie}
\end{equation}
Moreover, we assume that there exists a semimartingale $X$ with the form $%
X_{t}=X_{0}+\int_{0}^{t}\phi _{s}dB_{s}+V_{t}$, where $V$ is a process with
finite variation with $V=V^{+}-V^{-}$ and $V^{\pm }\in \mathcal{A}^{2}(0,T)$%
, satisfies for $t\in [0,T]$,
\[
L_{t}\leq X_{t}\leq U_{t}.
\]

\begin{definition}
\label{def1}We say a triple $(y_{t},z_{t},K_{t})_{0\leq t\leq T}\in \mathcal{%
S}^{2}(0,T)\times \mathcal{H}_{d}^{2}(0,T)\times \mathcal{A}^{2}(0,T)$ is a
solution of reflected BSDE with one barrier associated to $(\xi ,g,L)$, if
it satisfies

(i) for $t\in [0,T]$,
\[
y_{t}=\xi
+\int_{t}^{T}g(s,y_{s},z_{s})ds+K_{T}-K_{t}-\int_{t}^{T}z_{s}dB_{s};
\]

(ii) $y_{t}\geq L_{t}$, a.s. $0\leq t\leq T$;

(iii) $\int_{0}^{T}(y_{t}-L_{t})dK_{t}=0$.
\end{definition}

\begin{definition}
\label{def2}We say a quadruple $(y_{t},z_{t},A_{t},K_{t})_{0\leq t\leq T}$
is a solution of reflected BSDE with two barriers associated to $(\xi ,g,L,U)
$, for $y\in \mathcal{S}^{2}(0,T)$, $z\in \mathcal{H}_{d}^{2}(0,T)$ and $%
K^{\pm }$ $\in \mathcal{A}^{2}(0,T)$, if it satisfies

(i) for $t\in [0,T]$,
\[
y_{t}=\xi
+\int_{t}^{T}g(s,y_{s},z_{s})ds+K_{T}^{+}-K_{t}^{+}-(K_{T}^{-}-K_{t}^{-})-%
\int_{t}^{T}z_{s}dB_{s};
\]

(ii) $L_{t}\leq y_{t}\leq U_{t}$, a.s. $0\leq t\leq T$;

(iii) $\int_{0}^{T}(y_{t}-L_{t})dK_{t}^{+}=%
\int_{0}^{T}(y_{t}-U_{t})dK_{t}^{-}=0$.
\end{definition}

\section{Preliminaries}

In this section, we present some known results of reflected BSDE under
linear increasing condition. For $n\in \mathbf{N}$, we define
\begin{eqnarray}
\underline{g}_{n}(t,y,z) &=&\inf_{u\in \mathbf{R,v\in R}^{d}}\{g(t,u,v)+n(%
\left| y-u\right| +\left| z-v\right| )\},  \label{approximate} \\
\overline{g}_{n}(t,y,z) &=&\sup_{u\in \mathbf{R,v\in R}^{d}}\{g(t,u,v)-n(%
\left| y-u\right| +\left| z-v\right| )\}.  \nonumber
\end{eqnarray}
These sequences give important approximations of continuous functions by
Lipschitz functions (see proof in \cite{LS04}):

\begin{lemma}
\label{convegence}Let $g$ satisfy (H1), (H2) and (H3). Set $\mu =\max
\{\beta ,A\}$, then for $n>\mu $, we have for $t\in [0,T]$, $y\in \mathbf{R}$%
, $z\in \mathbf{R}^{d}$

(i) $-\mu (\left| y\right| +\left| z\right| +1)\leq \underline{g}%
_{n}(t,y,z)\leq g(t,y,z)\leq \overline{g}_{n}(t,y,z)\leq \mu (\left|
y\right| +\left| z\right| +1)$;

(ii) $\underline{g}_{n}(t,y,z)$ (resp. $\overline{g}_{n}(t,y,z)$) is
non-decreasing (resp. non-increasing) in $n$;

(iii) $\underline{g}_{n}(t,y,z)$ and $\overline{g}_{n}(t,y,z)$ are Lipschitz
in $(y,z)$ with parameter $n$, uniformly with respect to $(t,\omega )$;

(iv) If $(y_{n},z_{n})\rightarrow (y,z)$, as $n\rightarrow \infty $, then $%
\underline{g}_{n}(t,y_{n},z_{n})$ (resp. $\overline{g}_{n}(t,y_{n},z_{n}))$)
$\rightarrow g(t,y,z)$, as $n\rightarrow \infty $.
\end{lemma}

For $\xi \in \mathbf{L}^{2}(\mathcal{F}_{T})$ and $L$ satisfying (iii), we
consider two sequence of reflected BSDE associated to $(\xi ,\overline{g}%
_{n},L)$ and $(\xi ,\underline{g}_{n},L)$, for $n\in \mathbf{N}$,
respectively. For $n>\mu $, $\overline{g}_{n}$ and $\underline{g}_{n}$ are
Lipschitz in $(y,z)$, so reflected BSDEs have unique solutions in $\mathcal{S%
}^{2}(0,T)\times \mathcal{H}_{d}^{2}(0,T)\times \mathcal{A}^{2}(0,T)$,
denoted as $(\overline{y}^{n},\overline{z}^{n},\overline{K}^{n})$ and $(%
\underline{y}^{n},\underline{z}^{n},\underline{K}^{n})$ respectively.
Namely, the followings are satisfied:
\begin{eqnarray}
\overline{y}_{t}^{n} &=&\xi +\int_{t}^{T}\overline{g}_{n}(s,\overline{y}%
_{t}^{n},\overline{z}_{s}^{n})ds+\overline{K}_{T}^{n}-\overline{K}%
_{t}^{n}-\int_{t}^{T}\overline{z}_{s}^{n}dB_{s},  \label{s-above} \\
\overline{y}_{t}^{n} &\geq &L_{t},\;\;\int_{0}^{T}(\overline{y}%
_{t}^{n}-L_{t})d\overline{K}_{t}^{n}=0;  \nonumber
\end{eqnarray}
and
\begin{eqnarray}
\underline{y}_{t}^{n} &=&\xi +\int_{t}^{T}\underline{g}_{n}(s,\underline{y}%
_{t}^{n},\underline{z}_{s}^{n})ds+\underline{K}_{T}^{n}-\underline{K}%
_{t}^{n}-\int_{t}^{T}\underline{z}_{s}^{n}dB_{s},  \label{s-down} \\
\underline{y}_{t}^{n} &\geq &L_{t},0\leq t\leq T,\;\;\;\;\int_{0}^{T}(%
\underline{y}_{t}^{n}-L_{t})d\underline{K}_{t}^{n}=0.  \nonumber
\end{eqnarray}

We can prove

\begin{lemma}
\label{estimate}Under assumption (i)-(iii), there exists a constant $M_{0}$
independent of $n$, such that
\[
E[\sup_{0\leq t\leq T}\left| \underline{y}_{t}^{n}\right|
^{2}+\int_{0}^{T}\left| \underline{z}_{t}^{n}\right| ^{2}dt+\sup_{0\leq
t\leq T}\left| \overline{y}_{t}^{n}\right| ^{2}+\int_{0}^{T}\left| \overline{%
z}_{t}^{n}\right| ^{2}dt]\leq M_{0}.
\]
\end{lemma}

Then, we have the existence of the maximal and minimal solution.

\begin{theorem}
\label{exist}Let $(\xi ,g,L)$ be a triple satisfying the above assumptions,
in particular (i)-(iii). Then the reflected BSDE associated to $(\xi ,g,L)$,
has the maximal solution $(\overline{y}_{t},\overline{z}_{t},\overline{K}%
_{t})_{0\leq t\leq T}$ (resp. minimal solution $(\underline{y}_{t},%
\underline{z}_{t},\underline{K}_{t})_{0\leq t\leq T}$), i.e. it satisfies
definition \ref{def1}. Moreover $t\in [0,T]$,
\[
\underline{y}_{t}^{n}\leq \underline{y}_{t}^{n+1}\leq \underline{y}_{t}\leq
\overline{y}_{t}\leq \overline{y}_{t}^{n+1}\leq \overline{y}_{t}.
\]
And $(\overline{y}^{n},\overline{z}^{n},\overline{K}^{n})\rightarrow (%
\overline{y},\overline{z},\overline{K})$ and $(\underline{y}^{n},\underline{z%
}^{n},\underline{K}^{n})\rightarrow (\underline{y},\underline{z},\underline{K%
})$ both in $\mathcal{S}^{2}(0,T)\times \mathcal{H}_{d}^{2}(0,T)\times
\mathcal{S}^{2}(0,T)$ as $n\rightarrow \infty $.
\end{theorem}

\begin{remark}
In \cite{M97}, theorem 1 give the proof of the existence of the minimal
solution by the convergence of $(\underline{y}^{n},\underline{z}^{n},%
\underline{K}^{n})$. Due to the maximal solution, we can prove its existence
by considering $(\overline{y}^{n},\overline{z}^{n},\overline{K}^{n})$
symmetrically.
\end{remark}

\begin{remark}
We say the solution $(\overline{y}_{t},\overline{z}_{t},\overline{K}%
_{t})_{0\leq t\leq T}$ (resp. $(\underline{y}_{t},\underline{z}_{t},%
\underline{K}_{t})_{0\leq t\leq T}$) is maximal (resp. minimal) in the sense
that, if there exists another triple $(y_{t}^{\prime },z_{t}^{\prime
},K_{t}^{\prime })_{0\leq t\leq T}$ satisfies definition \ref{def}, then $%
y_{t}^{\prime }\leq \overline{y}_{t}$ (resp. $\underline{y}_{t}\leq
y_{t}^{\prime }$ ) $0\leq t\leq T$. This is an easy result of the following
comparison theorem.
\end{remark}

We have the following comparison result:

\begin{theorem}
\label{comparison}For $i=1,2$, assume $(\xi ^{i},g^{i},L^{i})$ satisfying
the assumptions (i)-(iii). Let $(\overline{y}_{t}^{i},\overline{z}_{t}^{i},%
\overline{K}_{t}^{i})_{0\leq t\leq T}$ be the maximal solution of reflected
BSDE associated to $(\xi ^{i},g^{i},L^{i})$. Moreover if for $(t,y,z)\in
[0,T]\times \mathbf{R}\times \mathbf{R}^{d}$, we have
\[
\xi ^{1}\geq \xi ^{2},g^{1}(t,y,z)\geq g^{2}(t,y,z),L_{t}^{1}\geq L_{t}^{2},
\]
then $\overline{y}_{t}^{1}\geq \overline{y}_{t}^{2}$, $t\in [0,T]$.
Furthermore, if $L^{1}=L^{2}$, then $\overline{K}_{t}^{1}-\overline{K}%
_{s}^{1}\leq \overline{K}_{t}^{2}-\overline{K}_{s}^{2}$, $0\leq s\leq t\leq T
$.
\end{theorem}

\begin{remark}
The comparison results still hold when we consider the minimal solutions $(%
\underline{y}_{t}^{i},\underline{z}_{t}^{i},\underline{K}_{t}^{i})_{0\leq
t\leq T}$ of corresponding reflected BSDEs.
\end{remark}

\begin{remark}
The proof of this theorem can be found in \cite{KLQT}, which follows from
the approximation in the proof of theorem \ref{exist} and the comparison
theorem for Lipschitz case in \cite{EKPPQ} and \cite{HLM}.
\end{remark}

Then we recall a continuous dependence result of reflected BSDE under
Lipschitz condition, which is proved in \cite{EKPPQ}.

\begin{theorem}
\label{depend}Let $(\xi ^{i},g^{i},L)$ be two triple satisfying the
assumptions (i), (iii) and $g^{i}$ is Lipschitz in $(y,z)$ uniformly in $%
(t,\omega )$, i.e. for some $k>0$, $t\in [0,T]$, $y$, $y^{\prime }\in
\mathbf{R}$, $z$, $z^{\prime }\in \mathbf{R}^{d}$, such that
\[
\left| g^{i}(t,y,z)-g^{i}(t,y^{\prime },z^{\prime })\right| \leq k(\left|
y-y^{\prime }\right| +\left| z-z^{\prime }\right| ).
\]
Suppose $(\overline{y}_{t}^{i},\overline{z}_{t}^{i},\overline{K}%
_{t}^{i})_{0\leq t\leq T}$ is the solution of reflected BSDE$(\xi
^{i},g^{i},L^{i})$, for $i=1,2$. Set
\begin{eqnarray*}
\bigtriangleup \xi  &=&\xi ^{1}-\xi ^{2},\bigtriangleup g=g^{1}-g^{2}, \\
\bigtriangleup y &=&y^{1}-y^{2},\bigtriangleup z=z^{1}-z^{2},\bigtriangleup
K=K^{1}-K^{2}.
\end{eqnarray*}
Then there exists a constant $C$ such that
\[
E[\sup_{0\leq t\leq T}\left| \bigtriangleup y_{t}\right|
^{2}+\int_{0}^{T}\left| \bigtriangleup z_{s}\right| ^{2}ds+\left|
\bigtriangleup K_{T}\right| ^{2}]\leq CE[\left| \bigtriangleup \xi \right|
^{2}+\int_{0}^{T}\left| \bigtriangleup g(s,y_{s}^{2},z_{s}^{2})\right| ds].
\]
\end{theorem}

For reflected BSDE with two continuous barriers, we have similar
results of existence of maximin and minimum solution and comparison
theorem.

\begin{theorem}
\label{exist2}Consider reflected BSDE with two barriers associated to $(\xi
,g,L,U)$, which satisfies assumptions (i)-(iii), it has the maximal solution
$(\overline{y}_{t},\overline{z}_{t},\overline{K}_{t}^{+},\overline{K}%
_{t}^{-})_{0\leq t\leq T}$ (resp. minimal solution $(\underline{y}_{t},%
\underline{z}_{t},\underline{K}_{t}^{+},\underline{K}_{t}^{-})_{0\leq t\leq
T}$), i.e. it satisfies definition \ref{def1}. Moreover $t\in [0,T]$,
\[
\underline{y}_{t}^{n}\leq \underline{y}_{t}^{n+1}\leq \underline{y}_{t}\leq
\overline{y}_{t}\leq \overline{y}_{t}^{n+1}\leq \overline{y}_{t}.
\]
And $(\overline{y}^{n},\overline{z}^{n},\overline{K}^{n+},\overline{K}%
^{n-})\rightarrow (\overline{y},\overline{z},\overline{K}^{+},\overline{K}%
^{-})$ and $(\underline{y}^{n},\underline{z}^{n},\underline{K}^{n+},%
\underline{K}^{n-})\rightarrow (\underline{y},\underline{z},\underline{K}%
^{+},\underline{K}^{-})$ both in $\mathcal{S}^{2}(0,T)\times \mathcal{H}%
_{d}^{2}(0,T)\times \mathcal{S}^{2}(0,T)$ as $n\rightarrow \infty $. Here $(%
\overline{y}^{n},\overline{z}^{n},\overline{K}^{n+},\overline{K}^{n-})$
(resp. $(\underline{y}^{n},\underline{z}^{n},\underline{K}^{n+},\underline{K}%
^{n-})$) is solution of reflected BSDE with two barriers associated to $(\xi
,\overline{g}_{n},L,U)$ (resp. $(\xi ,\underline{g}_{n},L,U)$ ).
\end{theorem}

\begin{theorem}
\label{comp2r}For $i=1,2$, assume $(\xi ^{i},g^{i},L^{i},U^{i})$ satisfying
the assumptions (i)-(iii). Let $(\overline{y}_{t}^{i},\overline{z}_{t}^{i},%
\overline{K}_{t}^{i+},\overline{K}_{t}^{i-})_{0\leq t\leq T}$ be the maximal
solution of reflected BSDE associated to $(\xi ^{i},g^{i},L^{i},U^{i})$.
Moreover if for $(t,y,z)\in [0,T]\times \mathbf{R}\times \mathbf{R}^{d}$, we
have
\[
\xi ^{1}\geq \xi ^{2},g^{1}(t,y,z)\geq g^{2}(t,y,z),L_{t}^{1}\geq
L_{t}^{2},U_{t}^{1}\geq U_{t}^{2},
\]
then $\overline{y}_{t}^{1}\geq \overline{y}_{t}^{2}$, $t\in [0,T]$.
Furthermore, if $L^{1}=L^{2}$, $U^{1}=U^{2}$, then $\overline{K}_{t}^{1+}-%
\overline{K}_{s}^{1+}\leq \overline{K}_{t}^{2+}-\overline{K}_{s}^{2+}$, $%
\overline{K}_{t}^{1-}-\overline{K}_{s}^{1-}\geq \overline{K}_{t}^{2-}-%
\overline{K}_{s}^{2-}$, $0\leq s\leq t\leq T$.
\end{theorem}

\begin{remark}
The comparison results still hold when we consider the minimal solutions $(%
\underline{y}_{t}^{i},\underline{z}_{t}^{i},\underline{K}_{t}^{i+},%
\underline{K}_{t}^{i-})_{0\leq t\leq T}$ of corresponding reflected BSDEs.
\end{remark}

\section{The solutions between $\overline{y}_{t}$ and $\protect\underline{y}%
_{t}$ for BSDE with one lower barrier}

Our first result is that between the maximal solution $\overline{y}_{t}$ and
the minimal solution $\underline{y}_{t}$, we can construct as many solution $%
(y_{t},z_{t},K_{t})_{0\leq t\leq T}$ as we want to satisfy the reflected BSDE%
$(\xi ,g,L)$, i.e. definition \ref{def1} is satisfied.

\begin{theorem}
\label{constuc}Assume that (i)-(iii) hold for $(\xi ,g,L)$. Let $(\underline{%
y}_{t},\underline{z}_{t},\underline{K}_{t})_{0\leq t\leq T}$ and $(\overline{%
y}_{t},\overline{z}_{t},\overline{K}_{t})_{0\leq t\leq T}$ be the minimal
and maximal solution of reflected BSDE$(\xi ,g,L)$, respectively, i.e.
definition \ref{def1} is satisfied for both triples. Then for any $t_{0}\in
[0,T]$, and $\eta \in \mathbf{L}^{2}(\mathcal{F}_{t_{0}})$, such that
\[
\underline{y}_{t_{0}}\leq \eta \leq
\overline{y}_{t_{0}},\;\;\;\mbox{a.s.,}
\]
there exists at least one solution $(y_{t},z_{t},K_{t})_{0\leq t\leq T}$ in $%
\mathcal{S}^{2}(0,T)\times \mathcal{H}_{d}^{2}(0,T)\times \mathcal{A}%
^{2}(0,T)$ of reflected BSDE$(\xi ,g,L)$ passing through $(t_{0},\eta )$,
namely $y_{t_{0}}=\eta $.
\end{theorem}

\noindent \textbf{Proof.}
On the interval $[0,t_{0}]$, since $\eta \in \mathbf{L}^{2}(\mathcal{F}%
_{t_{0}})$, there exists at least one triple $%
(y_{t}^{1},z_{t}^{1},K_{t}^{1})_{0\leq t\leq t_{0}}$ to be a solution of
reflected BSDE$(\eta ,g,L)$ on $[0,t_{0}]$, i.e.
\begin{eqnarray*}
y_{t}^{1} &=&\eta
+\int_{t}^{t_{0}}g(s,y_{s}^{1},z_{s}^{1})ds+K_{t_{0}}^{1}-K_{t}^{1}-%
\int_{t}^{t_{0}}z_{s}^{1}dB_{s}, \\
y_{t}^{1} &\geq &L_{t},0\leq t\leq
t_{0},\;\;\;\;\int_{0}^{t_{0}}(y_{t}^{1}-L_{t})dK_{t}^{1}=0.
\end{eqnarray*}
Then we fix a process $z^{2}\in \mathcal{H}_{d}^{2}(t_{0},T)$, and consider
a (strong) solution $(y_{t}^{2})_{t_{0}\leq t\leq T}$ of the following SDE
\[
y_{t}^{2}=\eta
-\int_{t_{0}}^{t}g(s,y_{s}^{2},z_{s}^{2})ds+\int_{t_{0}}^{t}z_{s}^{2}dB_{s}.
\]

Define a stopping time $\tau =\inf \{t\geq t_{0},y_{t}^{2}\notin (\underline{%
y}_{t},\overline{y}_{t})\}$. Since $\underline{y}_{T}=\overline{y}_{T}$, we
get $\tau <T$. Notice that $y_{t_{0}}^{2}=\eta \geq L_{t_{0}}$, and by the
continuity of solutions we know that on the interval $[t_{0},\tau ]$, $%
y_{t}^{2}\geq \underline{y}_{t}\geq L_{t}$, which implies that on this
interval $(y_{t}^{2},z_{t}^{2},0)_{t_{0}\leq t\leq \tau }$ is also a
solution to reflected BSDE$(\xi ,g,L)$ on $[t_{0},\tau ]$.

Now we denote the triple on $[0,T]$%
\begin{eqnarray*}
y_{t} &=&y_{t}^{1}1_{[0,t_{0})}(t)+y_{t}^{2}1_{[t_{0},\tau )}(t)+\overline{y}%
_{t}1_{[\tau ,T]}(t)1_{\{y_{\tau }=\overline{y}_{\tau }\}}+\underline{y}%
_{t}1_{[\tau ,T]}(t)1_{\{y_{\tau }<\overline{y}_{\tau }\}}, \\
z_{t} &=&z_{t}^{1}1_{[0,t_{0})}(t)+z_{t}^{2}1_{[t_{0},\tau )}(t)+\overline{z}%
_{t}1_{[\tau ,T]}(t)1_{\{y_{\tau }=\overline{y}_{\tau }\}}+\underline{z}%
_{t}1_{[\tau ,T]}(t)1_{\{y_{\tau }<\overline{y}_{\tau }\}}, \\
K_{t} &=&K_{t}^{1}1_{[0,t_{0})}(t)+K_{t_{0}}^{1}1_{[t_{0},\tau )}(t)+%
\overline{K}_{t}1_{[\tau ,T]}(t)1_{\{y_{\tau }=\overline{y}_{\tau }\}}+%
\underline{K}_{t}1_{[\tau ,T]}(t)1_{\{y_{\tau }<\overline{y}_{\tau }\}}.
\end{eqnarray*}
It is easy to check that $(y_{t},z_{t},K_{t})_{0\leq t\leq T}\in \mathcal{S}%
^{2}(0,T)\times \mathcal{H}_{d}^{2}(0,T)\times \mathcal{A}^{2}(0,T)$
satisfies definition \ref{def1}, which means that it is a solution of
reflected BSDE associated to $(\xi ,g,L)$. $\square $

\begin{remark}
This result is still true if we replace the linear increasing property of $g$
(i.e. (H2)) by quadratic growth assumption of \cite{KLQT}.
\end{remark}

\section{Similar result for reflected BSDE with two barriers}

We have similar result for reflcted BSDE with two barriers for constructing
as many solution $(y_{t},z_{t},K_{t}^{+},K_{t}^{-})_{0\leq t\leq T}$ as we
want to satisfy the reflected BSDE$(\xi ,g,L,U)$, i.e. definition \ref{def2}
is satisfied.

\begin{theorem}
Assume that (i)-(iii) hold for $(\xi ,g,L,U)$. Let $(\underline{y}_{t},%
\underline{z}_{t},\underline{K}_{t}^{+},\underline{K}_{t}^{-})_{0\leq t\leq
T}$ and $(\overline{y}_{t},\overline{z}_{t},\overline{K}_{t}^{+},\overline{K}%
_{t}^{-})_{0\leq t\leq T}$ be the minimal and maximal solution of reflected
BSDE$(\xi ,g,L,U)$, respectively, i.e. definition \ref{def2} is satisfied
for both triples. Then for any $t_{0}\in [0,T]$, and $\eta \in \mathbf{L}%
^{2}(\mathcal{F}_{t_{0}})$, such that
\[
\underline{y}_{t_{0}}\leq \eta \leq
\overline{y}_{t_{0}},\;\;\;\mbox{a.s.,}
\]
there exists at least one solution $(y_{t},z_{t},K_{t}^{+},K_{t}^{-})_{0\leq
t\leq T}$ in $\mathcal{S}^{2}(0,T)\times \mathcal{H}_{d}^{2}(0,T)\times (%
\mathcal{A}^{2}(0,T))^{2}$ of reflected BSDE$(\xi ,g,L,U)$ passing through $%
(t_{0},\eta )$, namely $y_{t_{0}}=\eta $.
\end{theorem}

\noindent \textbf{Proof.}
On the interval $[0,t_{0}]$, since $\eta \in \mathbf{L}^{2}(\mathcal{F}%
_{t_{0}})$, there exists at least one triple $%
(y_{t}^{1},z_{t}^{1},K_{t}^{1+},K_{t}^{1-})_{0\leq t\leq t_{0}}$ to be a
solution of reflected BSDE$(\eta ,g,L,U)$ on $[0,t_{0}]$, i.e.
\begin{eqnarray*}
y_{t}^{1} &=&\eta
+\int_{t}^{t_{0}}g(s,y_{s}^{1},z_{s}^{1})ds+K_{t_{0}}^{1+}-K_{t}^{1+}-(K_{t_{0}}^{1-}-K_{t}^{1-})-\int_{t}^{t_{0}}z_{s}^{1}dB_{s},
\\
L_{t} &\leq &y_{t}^{1}\leq U_{t},0\leq t\leq
t_{0},\;\;\;\;\int_{0}^{t_{0}}(y_{t}^{1}-L_{t})dK_{t}^{1+}=%
\int_{0}^{t_{0}}(y_{t}^{1}-U_{t})dK_{t}^{1-}=0.
\end{eqnarray*}
Then we fix a process $z^{2}\in \mathcal{H}_{d}^{2}(t_{0},T)$, and consider
a (strong) solution $(y_{t}^{2})_{t_{0}\leq t\leq T}$ of the following SDE
\[
y_{t}^{2}=\eta
-\int_{t_{0}}^{t}g(s,y_{s}^{2},z_{s}^{2})ds+\int_{t_{0}}^{t}z_{s}^{2}dB_{s}.
\]

Define a stopping time $\tau =\inf \{t\geq t_{0},y_{t}^{2}\notin (\underline{%
y}_{t},\overline{y}_{t})\}$. Since $\underline{y}_{T}=\overline{y}_{T}$, we
get $\tau <T$. Notice that $L_{t_{0}}\leq y_{t_{0}}^{2}=\eta \leq U_{t_{0}}$%
, and by the continuity of solutions we know that on the interval $%
[t_{0},\tau ]$, $L_{t}\leq \underline{y}_{t}\leq y_{t}^{2}\leq \overline{y}%
_{t}\leq U_{t}$, which implies that on this interval $%
(y_{t}^{2},z_{t}^{2},0,0)_{t_{0}\leq t\leq \tau }$ is also a solution to
reflected BSDE$(\xi ,g,L,U)$ on $[t_{0},\tau ]$.

Now we denote the triple on $[0,T]$%
\begin{eqnarray*}
y_{t} &=&y_{t}^{1}1_{[0,t_{0})}(t)+y_{t}^{2}1_{[t_{0},\tau )}(t)+\overline{y}%
_{t}1_{[\tau ,T]}(t)1_{\{y_{\tau }=\overline{y}_{\tau }\}}+\underline{y}%
_{t}1_{[\tau ,T]}(t)1_{\{y_{\tau }<\overline{y}_{\tau }\}}, \\
z_{t} &=&z_{t}^{1}1_{[0,t_{0})}(t)+z_{t}^{2}1_{[t_{0},\tau )}(t)+\overline{z}%
_{t}1_{[\tau ,T]}(t)1_{\{y_{\tau }=\overline{y}_{\tau }\}}+\underline{z}%
_{t}1_{[\tau ,T]}(t)1_{\{y_{\tau }<\overline{y}_{\tau }\}}, \\
K_{t}^{+} &=&K_{t}^{1+}1_{[0,t_{0})}(t)+K_{t_{0}}^{1+}1_{[t_{0},\tau )}(t)+%
\overline{K}_{t}^{+}1_{[\tau ,T]}(t)1_{\{y_{\tau }=\overline{y}_{\tau }\}}+%
\underline{K}_{t}^{+}1_{[\tau ,T]}(t)1_{\{y_{\tau }<\overline{y}_{\tau }\}},
\\
K_{t}^{-} &=&K_{t}^{1-}1_{[0,t_{0})}(t)+K_{t_{0}}^{1-}1_{[t_{0},\tau )}(t)+%
\overline{K}_{t}^{-}1_{[\tau ,T]}(t)1_{\{y_{\tau }=\overline{y}_{\tau }\}}+%
\underline{K}_{t}^{-}1_{[\tau ,T]}(t)1_{\{y_{\tau }<\overline{y}_{\tau }\}}
\end{eqnarray*}
It is easy to check that $(y_{t},z_{t},K_{t}^{+},K_{t}^{-})_{0\leq t\leq
T}\in \mathcal{S}^{2}(0,T)\times \mathcal{H}_{d}^{2}(0,T)\times (\mathcal{A}%
^{2}(0,T))^{2}$ satisfies definition \ref{def2}, which means that it is a
solution of reflected BSDE associated to $(\xi ,g,L,U)$. $\square $

\begin{remark}
This result is still true if we replace the linear increasing
property of $g$ (i.e. (H2)) by quadratic growth assumption as in
\cite{KLQT}.
\end{remark}

\section{One uniqueness theorem when $g$ only depend on $z$}

In this section, we will prove one uniqueness result of the reflected BSDE
whose coefficient is only depend on $z$. We still assume (H1), (H2) and (H3)
hold for $g:[0,T]\times \Omega \mathbf{\times R}^{d}\rightarrow \mathbf{R}$,
however in order to get the uniqueness of the solution, we need

(H4) uniform continuity: $g(t,\cdot )$ is uniformly continuous in $z$,
uniformly with respect to $(\omega ,t)$. More precisely, there exists a
continuous, non-decreasing function $\phi :\mathbf{R}^{+}\rightarrow \mathbf{%
R}^{+}$ with linear growth with parameter $A$ and satisfying $\phi (0)=0$
such that for any $t\in [0,T]$, $z_{1}$, $z_{2}\in \mathbf{R}^{d}$,
\[
\left| g(t,z_{1})-g(t,z_{2})\right| \leq \phi (\left| z_{1}-z_{2}\right| ).
\]
In fact, this assumption implies assumption (H3).

Then as (\ref{approximate}) we define
\begin{eqnarray*}
\underline{g}_{n}(z) &=&\inf_{u\in \mathbf{R}^{d}}\{g(u)+n\left| z-u\right|
\}, \\
\overline{g}_{n}(z) &=&\sup_{u\in \mathbf{R}^{d}}\{g(u)-n\left| z-u\right|
\}.
\end{eqnarray*}
So for $n>\mu $, Lemma \ref{convegence} also holds for $\underline{g}_{n}(z)$
and $\overline{g}_{n}(z)$. Furthermore we have

\begin{lemma}
\label{g-control}For $z\in \mathbf{R}^{d}$, $0\leq g(z)-\underline{g}%
_{n}(z)\leq \phi (\frac{\mu }{n-\mu })$ and $0\leq \overline{g}%
_{n}(z)-g(z)\leq \phi (\frac{\mu }{n-\mu })$.
\end{lemma}

This lemma is proved in \cite{G06}.

We first consider reflected BSDE with one barrier. For $\xi \in \mathbf{L}%
^{2}(\mathcal{F}_{T})$ and $L$ satisfying (iii), we consider two sequence of
reflected BSDE associated to $(\xi ,\overline{g}_{n},L)$ and $(\xi ,%
\underline{g}_{n},L)$, for $n\in \mathbf{N}$, respectively. For $n>\mu $, $%
\overline{g}_{n}$ and $\underline{g}_{n}$ are Lipschitz in $z$, so reflected
BSDEs have unique solutions, which are denoted as $(\overline{y}^{n},%
\overline{z}^{n},\overline{K}^{n})$ and $(\underline{y}^{n},\underline{z}%
^{n},\underline{K}^{n})$ respectively. So the followings are satisfied:
\begin{eqnarray}
\overline{y}_{t}^{n} &=&\xi +\int_{t}^{T}\overline{g}_{n}(s,\overline{z}%
_{s}^{n})ds+\overline{K}_{T}^{n}-\overline{K}_{t}^{n}-\int_{t}^{T}\overline{z%
}_{s}^{n}dB_{s},  \label{sup-s} \\
\overline{y}_{t}^{n} &\geq &L_{t},\;\;\int_{0}^{T}(\overline{y}%
_{t}^{n}-L_{t})d\overline{K}_{t}^{n}=0;  \nonumber
\end{eqnarray}
and
\begin{eqnarray}
\underline{y}_{t}^{n} &=&\xi +\int_{t}^{T}\underline{g}_{n}(s,\underline{z}%
_{s}^{n})ds+\underline{K}_{T}^{n}-\underline{K}_{t}^{n}-\int_{t}^{T}%
\underline{z}_{s}^{n}dB_{s},  \label{sub-s} \\
\underline{y}_{t}^{n} &\geq &L_{t},0\leq t\leq T,\;\;\;\;\int_{0}^{T}(%
\underline{y}_{t}^{n}-L_{t})d\underline{K}_{t}^{n}=0.  \nonumber
\end{eqnarray}

Similarly to \cite{G06}, we have the following lemma.

\begin{lemma}
\label{diff}Under assumption (H1), (H2) and (H4), let $(\underline{y}_{t},%
\underline{z}_{t},\underline{K}_{t})_{0\leq t\leq T}$ (resp. $(\overline{y}%
_{t},\overline{z}_{t},\overline{K}_{t})_{0\leq t\leq T}$) be the minimal
(resp. maximal) solution of reflected BSDE$(\xi ,g,L)$, i.e. they satisfy
definition \ref{def1}. Then we have for $n>\mu $, $E[\overline{y}_{t}^{n}-%
\underline{y}_{t}^{n}]\leq 2\phi (\frac{\mu }{n-\mu })T$.
\end{lemma}

\noindent \textbf{Proof.} From (\ref{sup-s}) and (\ref{sub-s}), we
get
\begin{eqnarray*}
\overline{y}_{t}^{n}-\underline{y}_{t}^{n} &=&\int_{t}^{T}(\overline{g}%
_{n}(s,\overline{z}_{s}^{n})-\underline{g}_{n}(s,\underline{z}%
_{s}^{n}))ds-\int_{t}^{T}(\overline{z}_{s}^{n}-\underline{z}_{s}^{n})dB_{s}
\\
&&+(\overline{K}_{T}^{n}-\overline{K}_{t}^{n})-(\underline{K}_{T}^{n}-%
\underline{K}_{t}^{n}).
\end{eqnarray*}
Since $\overline{g}_{n}(t,z)\geq \underline{g}_{n}(t,z)$, for $(t,z)\in
[0,T]\times \mathbf{R}^{d}$, by the comparison theorem in \cite{HLM} or
Theorem \ref{comparison}, we have
\[
0\leq \overline{K}_{T}^{n}-\overline{K}_{t}^{n}\leq \underline{K}_{T}^{n}-%
\underline{K}_{t}^{n}.
\]
Then
\begin{equation}
0\leq \overline{y}_{t}^{n}-\underline{y}_{t}^{n}\leq \widetilde{y}_{t}^{n},
\label{control-yn}
\end{equation}
where
\begin{equation}
\widetilde{y}_{t}^{n}:=\int_{t}^{T}(\overline{g}_{n}(s,\overline{z}_{s}^{n})-%
\underline{g}_{n}(s,\underline{z}_{s}^{n}))ds-\int_{t}^{T}(\overline{z}%
_{s}^{n}-\underline{z}_{s}^{n})dB_{s}.  \label{diff-equ}
\end{equation}

Set $\widetilde{z}_{t}^{n}:=\overline{z}_{s}^{n}-\underline{z}_{s}^{n}$ and
\begin{equation}
l_{t}^{n,j}=\left\{
\begin{array}{ll}
\frac{\underline{g}_{n}(t,\widehat{z}_{t}^{n,j-1})-\underline{g}_{n}(t,%
\widehat{z}_{t}^{n,j})}{\overline{z}_{t}^{n,j}-\underline{z}_{t}^{n,j}}, & %
\mbox{if}\hskip 0.1cm\overline{z}_{t}^{n,j}\neq \underline{z}_{t}^{n,j} \\
0, & \mbox{if}\hskip 0.1cm\overline{z}_{t}^{n,j}\neq \underline{z}_{t}^{n,j}
\end{array}
,\right.  \label{linear}
\end{equation}
where $\widehat{z}_{t}^{n,j}$ is the vector whose first $j$ components are
equal to those of $\underline{z}_{t}^{n}$ and whose of last $d-j$ components
are equal to those of $\overline{z}_{t}^{n}$, i.e. $\widehat{z}_{t}^{j}=(%
\underline{z}_{t}^{n,1},\ldots ,\underline{z}_{t}^{n,j},\overline{z}%
_{t}^{n,j+1},\ldots ,\overline{z}_{t}^{n,d})$. Here $l_{t}^{j}$ is $j$th
vector of $l$. Note that
\begin{eqnarray*}
\overline{g}_{n}(s,\overline{z}_{s}^{n})-\underline{g}_{n}(s,\underline{z}%
_{s}^{n}) &=&\overline{g}_{n}(s,\overline{z}_{s}^{n})-\underline{g}_{n}(s,%
\overline{z}_{s}^{n})+\underline{g}_{n}(s,\overline{z}_{s}^{n})-\underline{g}%
_{n}(s,\underline{z}_{s}^{n}) \\
&=&\widetilde{g}_{s}^{n}+\underline{g}_{n}(s,\overline{z}_{s}^{n})-%
\underline{g}_{n}(s,\underline{z}_{s}^{n}),
\end{eqnarray*}
where $\widetilde{g}_{s}^{n}:=\overline{g}_{n}(s,\overline{z}_{s}^{n})-%
\underline{g}_{n}(s,\overline{z}_{s}^{n})$, we know that $(\widetilde{y}%
_{t}^{n},\widetilde{z}_{t}^{n})$ satisfies the following linear BSDE on $%
[0,T]$,
\[
\widetilde{y}_{t}^{n}=\int_{t}^{T}(l_{s}^{n}\widetilde{z}_{s}^{n}+\widetilde{%
g}_{s}^{n})ds-\int_{t}^{T}\widetilde{z}_{s}^{n}dB_{s}.
\]
Since now $\underline{g}_{n}$ is Lipschitz in $z$ with parameter $n$, $%
\left| l_{s}^{n}\right| \leq n$. Consider the solution of linear SDE $%
q_{t}^{n}:=\exp [\int_{0}^{t}l_{s}^{n}dB_{s}-\frac{1}{2}\int_{0}^{t}\left|
l_{s}^{n}\right| ^{2}ds]$, applying It\^{o}'s formula to $q_{t}^{n}%
\widetilde{y}_{t}^{n}$ on $[t,T]$ and taking conditional expectation, then
we get
\begin{eqnarray*}
\widetilde{y}_{t}^{n} &=&(q_{t}^{n})^{-1}E[\int_{t}^{T}q_{s}^{n}\widetilde{g}%
_{s}^{n}ds|\mathcal{F}_{t}] \\
&=&E[\int_{t}^{T}(\int_{t}^{s}l_{s}^{n}dB_{s}-\frac{1}{2}\int_{t}^{s}\left|
l_{s}^{n}\right| ^{2}ds)\widetilde{g}_{s}^{n}ds|\mathcal{F}_{t}].
\end{eqnarray*}
By lemma \ref{g-control}, we have $0\leq \widetilde{g}_{s}^{n}\leq 2\phi (%
\frac{\mu }{n-\mu })$, for $s\in [0,T]$, $n>\mu $. Therefore

\begin{eqnarray}
E[\widetilde{y}_{t}^{n}] &=&E[\int_{t}^{T}(\int_{t}^{s}l_{s}^{n}dB_{s}-\frac{%
1}{2}\int_{t}^{s}\left| l_{s}^{n}\right| ^{2}ds)\widetilde{g}_{s}^{n}ds]
\label{control-e} \\
&\leq &2\phi (\frac{\mu }{n-\mu })E[\int_{t}^{T}(\int_{t}^{s}l_{s}^{n}dB_{s}-%
\frac{1}{2}\int_{t}^{s}\left| l_{s}^{n}\right| ^{2}ds)ds]  \nonumber \\
&\leq &2\phi (\frac{\mu }{n-\mu })T,  \nonumber
\end{eqnarray}
in view of $E[\int_{t}^{s}l_{s}^{n}dB_{s}-\frac{1}{2}\int_{t}^{s}\left|
l_{s}^{n}\right| ^{2}ds]=1$, for $t\leq s\leq T$, which follows from the
fact that $q^{n}$ is a exponential martingale. The result follows from (\ref
{control-yn}) and (\ref{control-e}). $\square $

With these preparations, we present our main result of this section.

\begin{theorem}
\label{unique1}Assume assumptions (i), (iii) hold for $\xi $ and $L$, and $%
g:[0,T]\times \mathbf{R}^{d}\rightarrow \mathbf{R}$ satisfies (H1), (H2) and
(H4). The solution of reflected BSDE$(\xi ,g,L)$ is unique.
\end{theorem}

\noindent \textbf{Proof.}
From Lemma \ref{diff}, we have $E[\overline{y}_{t}^{n}-\underline{y}%
_{t}^{n}]\leq 2\phi (\frac{\mu }{n-\mu })T$, for $n>\mu $. So $E[\overline{y}%
_{t}^{n}-\underline{y}_{t}^{n}]\rightarrow 0$, as $n\rightarrow \infty $.
While Lemma \ref{estimate} implies $(\overline{y}_{t}^{n}-\underline{y}%
_{t}^{n})$ is bounded in $\mathcal{S}^{2}(0,T)$ uniformly in $n$, we get $E[(%
\overline{y}_{t}^{n}-\underline{y}_{t}^{n})^{2}]\rightarrow 0$, as $%
n\rightarrow \infty $, in view of $\overline{y}_{t}^{n}-\underline{y}%
_{t}^{n}\geq 0$.

Let $(\underline{y}_{t},\underline{z}_{t},\underline{K}_{t})_{0\leq t\leq T}$
(resp. $(\overline{y}_{t},\overline{z}_{t},\overline{K}_{t})_{0\leq t\leq T}$%
) be the minimal (resp. maximal) solution of reflected BSDE$(\xi ,g,L)$, by
the convergence result of Theorem \ref{exist}, we obtain
\[
E[(\overline{y}_{t}-\underline{y}_{t})^{2}]\leq E[(\overline{y}_{t}-%
\overline{y}_{t}^{n})^{2}]+E[(\overline{y}_{t}^{n}-\underline{y}%
_{t}^{n})^{2}]+E[(\underline{y}_{t}^{n}-\underline{y}_{t})^{2}]\rightarrow 0,
\]
as $n\rightarrow \infty $. The proof is complete. $\square $

For reflected BSDE with two barriers, we have similar result.

\begin{theorem}
Assume assumptions (i), (iii) hold for $\xi $, $L$ and $U$ and $%
g:[0,T]\times \mathbf{R}^{d}\rightarrow \mathbf{R}$ satisfies (H1), (H2) and
(H4). The solution of reflected BSDE$(\xi ,g,L,U)$ is unique.
\end{theorem}

\noindent \textbf{Proof.} The proof is similar to the one of Theorem
\ref{unique1}. First we need to
prove a result as Lemma \ref{diff}: for $n>\mu $, $E[\overline{y}_{t}^{n}-%
\underline{y}_{t}^{n}]\leq 2\phi (\frac{\mu }{n-\mu })T$, where $\overline{y}%
_{t}^{n}$ (resp. $\underline{y}_{t}^{n}$) is the solution of reflected BSDE$%
(\xi ,\overline{g}^{n},L,U)$ (resp. $(\xi ,\underline{g}^{n},L,U)$). In
fact, when we consider the difference of $\overline{y}_{t}^{n}$ and $%
\underline{y}_{t}^{n}$.
\begin{eqnarray*}
\overline{y}_{t}^{n}-\underline{y}_{t}^{n} &=&\int_{t}^{T}(\overline{g}%
_{n}(s,\overline{z}_{s}^{n})-\underline{g}_{n}(s,\underline{z}%
_{s}^{n}))ds-\int_{t}^{T}(\overline{z}_{s}^{n}-\underline{z}_{s}^{n})dB_{s}
\\
&&+(\overline{K}_{T}^{n+}-\overline{K}_{t}^{n+})-(\underline{K}_{T}^{n+}-%
\underline{K}_{t}^{n+})-(\overline{K}_{T}^{n-}-\overline{K}_{t}^{n-})+(%
\underline{K}_{T}^{n-}-\underline{K}_{t}^{n-}).
\end{eqnarray*}
Thanks to Theorem \ref{comp2r}, we have $0\leq \overline{K}_{T}^{n+}-%
\overline{K}_{t}^{n+}\leq \underline{K}_{T}^{n+}-\underline{K}_{t}^{n+}$, $%
\overline{K}_{T}^{n-}-\overline{K}_{t}^{n-}\geq \underline{K}_{T}^{n-}-%
\underline{K}_{t}^{n-}\geq 0$. Then we get $0\leq \overline{y}_{t}^{n}-%
\underline{y}_{t}^{n}\leq \widetilde{y}_{t}^{n}$, where $\widetilde{y}%
_{t}^{n}$ is the solution of (\ref{diff-equ}). So following same method, it
follows  for $n>\mu $, $E[\overline{y}_{t}^{n}-\underline{y}_{t}^{n}]\leq
2\phi (\frac{\mu }{n-\mu })T$.

Let $(\underline{y}_{t},\underline{z}_{t},\underline{K}_{t}^{+},\underline{K}%
_{t}^{-})_{0\leq t\leq T}$ (resp. $(\overline{y}_{t},\overline{z}_{t},%
\overline{K}_{t}^{+},\overline{K}_{t}^{-})_{0\leq t\leq T}$) be the minimal
(resp. maximal) solution of reflected BSDE$(\xi ,g,L,U)$, by the convergence
result of Theorem \ref{exist2}, we obtain
\[
E[(\overline{y}_{t}-\underline{y}_{t})^{2}]\leq E[(\overline{y}_{t}-%
\overline{y}_{t}^{n})^{2}]+E[(\overline{y}_{t}^{n}-\underline{y}%
_{t}^{n})^{2}]+E[(\underline{y}_{t}^{n}-\underline{y}_{t})^{2}]\rightarrow 0,
\]
as $n\rightarrow \infty $. The proof is complete. $\square $

\section{Uniqueness results about the Disturbance of Coefficient}

In this section, we consider the uniqueness problem of solution for
reflected BSDEs with respect to a disturbance $c\in \mathbf{R}$ of its
coefficient as $g(t,y,z)+c$, under the uniform continuity assumption of $g$
on $(y,z)$. More precisely, we assume that the function $g:[0,T]\times
\Omega \times \mathbf{R\times R}^{d}\rightarrow \mathbf{R}$, satisfies

(H4') uniform continuity: $g(t,\cdot ,\cdot )$ is uniformly continuous in $%
(y,z)$, uniformly with respect to $(\omega ,t)$, i.e., there exists a
continuous, non-decreasing function $\phi :\mathbf{R}^{+}\rightarrow \mathbf{%
R}^{+}$ satisfying $\phi (0)=0$ and $\phi (x)\leq A(1+\left| x\right| )$
such that for any $t\in [0,T]$, $y_{1}$, $y_{2}\in \mathbf{R}$, $z_{1}$, $%
z_{2}\in \mathbf{R}^{d}$,
\[
\left| g(t,y_{1},z_{1})-g(t,y_{2},z_{2})\right| \leq \phi (\left|
y_{1}-y_{2}\right| +\left| z_{1}-z_{2}\right| ).
\]
In fact, this assumption implies assumption (H3).

We define $\underline{g}_{n}(t,y,z)$ and $\overline{g}_{n}(t,y,z)$ same as
in (\ref{approximate}), then we know that lemma \ref{convegence} holds for $%
\underline{g}_{n}$ and $\overline{g}_{n}$. Moreover we have the following
lemma, which is proved in \cite{G06-2}:

\begin{lemma}
\label{diff2}Let $g$ satisfy (H1), (H2) and (H4'). For $n\in \mathbf{N}$,
Set $\mu =\max \{\beta ,A\}$, then for $n>\mu $, we have for $t\in [0,T]$, $%
y\in \mathbf{R}$, $z\in \mathbf{R}^{d}$, $0\leq g(t,y,z)-\underline{g}%
_{n}(t,y,z)\leq \phi (\frac{\mu }{n-\mu })$ and $0\leq \overline{g}%
_{n}(t,y,z)-g(t,y,z)\leq \phi (\frac{\mu }{n-\mu })$.
\end{lemma}

Set $(\underline{y}_{t}^{c},\underline{z}_{t}^{c},\underline{K}%
_{t}^{c})_{0\leq t\leq T}$ (resp. $(\overline{y}_{t}^{c},\overline{z}%
_{t}^{c},\overline{K}_{t}^{c})_{0\leq t\leq T}$) be the minimal (resp.
maximal) solution of reflected BSDE$(\xi ,g+c,L)$. Before considering the
main result of this section, we prove a comparison theorem.

\begin{lemma}
\label{comp}Assume (H1), (H2) and (H4') hold for $g$. For a given $c>0$, $%
\xi \in \mathbf{L}^{2}(\mathcal{F}_{T})$ and $L$ satisfying (iii), set $%
(y_{t},z_{t},K_{t})_{0\leq t\leq T}$ be a solution of reflected BSDE$(\xi
,g,L)$, then we have
\[
\underline{y}_{t}^{c}\geq y_{t}\mbox{, }t\in [0,T].
\]
\end{lemma}

\noindent \textbf{Proof.}
By lemma \ref{diff2}, (i) of lemma \ref{convegence} and the continuity of $%
\phi $ at $x=0$, there exists an enough large $n_{0}>\mu $, such that $\phi (%
\frac{\mu }{n_{0}-\mu })<\frac{c}{2}$ and
\[
g(t,y,z)\leq \underline{g}_{n_{0}}(t,y,z)+\phi (\frac{\mu }{n_{0}-\mu })\leq
\underline{g}_{n_{0}}(t,y,z)+\frac{c}{2}<g(t,y,z)+c.
\]
Notice that $\underline{g}_{n_{0}}$ is Lipschitz in $(y,z)$ with parameter $%
n_{0}$, so the reflected BSDE$(\xi ,\underline{g}_{n_{0}}+\frac{c}{2},L)$
admits the unique solution defined as $%
(y_{t}^{c,n_{0}},z_{t}^{c,n_{0}},K_{t}^{c,n_{0}})$. By the comparison
theorem \ref{comparison}, we get
\[
y_{t}\leq \overline{y}_{t}\leq y_{t}^{c,n_{0}}\leq \underline{y}_{t}^{c}%
\mbox{, }t\in [0,T].
\]
$\square $

Now we introduce two auxiliary functions
\[
\overline{m}(t,c)=E[\overline{y}_{t}^{c}],\underline{m}(t,c)=E[\underline{y}%
_{t}^{c}],
\]
which have following properties:

\begin{proposition}
\label{prop-m}(i)For any $c\in \mathbf{R}$, $t\rightarrow \overline{m}(t,c)$
or $\underline{m}(t,c)$ is continuous;

(ii) for any $t\in [0,T]$, $c\rightarrow \overline{m}(t,c)$ or $\underline{m}%
(t,c)$ is nondecreasing;

(iii) for any $t\in [0,T]$, $c\rightarrow \overline{m}(t,c)$ is left
continuous and $c\rightarrow \underline{m}(t,c)$ is right continuous.
\end{proposition}

The proof is similar to the proposition 7 in \cite{G06-2}, with the helps of
comparison theorem \ref{comparison}, theorem \ref{depend} and the
convergence results of theorem \ref{exist}. So we omit here.

From the properties of $\overline{m}(t,c)$ and $\underline{m}(t,c)$ and
lemma \ref{comp}, we obtain a necessary and sufficient condition for the
uniqueness of solution of reflected BSDE under uniformly continuous property.

\begin{theorem}
\label{equivalent}Let $(\xi ,g,L)$ satisfy (i)-(iii) and (H4') hold for $g$.
For $c_{0}\in \mathbf{R}$, the following statements are equivalent:

(i) The reflected BSDE$(\xi ,g+c_{0},L)$ admits the unique solution;

(ii) $\underline{m}(t,c)$ is continuous at $c=c_{0}$, for all $t\in [0,T]$;

(iii) $\overline{m}(t,c)$ is continuous at $c=c_{0}$, for all $t\in [0,T]$;

(iv) $\overline{m}(t,c_{0})=\overline{m}(t,c_{0})$, for all $t\in [0,T]$.
\end{theorem}

\noindent \textbf{Proof.} It is easy to check that (i) and (iv) are
equivalent. Thanks to proposition \ref{prop-m} and comparison
theorem \ref{comparison}, with the similar proof in theorem 9 in
\cite{G06-2}, we deduce that (i)$\Rightarrow $(ii) or (iii).
Similarly, by lemma \ref{comp}, we can prove (ii) or
(iii)$\Rightarrow $(i). $\square $

Finally, we give the last result of this section.

\begin{theorem}
Let $(\xi ,g,L)$ satisfy (i)-(iii) and (H4') hold for $g$. Then the set of
real number $c$ such that the reflected BSDE with a disturbance $c\in
\mathbf{R}$ of its coefficient, i.e. reflected BSDE$(\xi ,g+c,L)$, admits
more than one solution, is at most countable.
\end{theorem}

\noindent \textbf{Proof.}
From theorem \ref{equivalent}, we deduce that $\overline{m}(t,c_{0})=%
\overline{m}(t,c_{0})$, for all $t\in [0,T]$, is equivalent to the
uniqueness of solution of reflected BSDE$(\xi ,g+c,L)$. Since $c\rightarrow
\overline{m}(t,c)$ or $\underline{m}(t,c)$ is monotone, it has at most
countable discontinuous points. While $t\rightarrow \overline{m}(t,c)$ or $%
\underline{m}(t,c)$ is continuous, our results follow by classical
techniques. $\square $

\end{document}